\documentclass{article}
\usepackage{amsmath,amssymb}

\begin{document}

\title{A long cycle theorem involving Fan-type degree condition and neighborhood intersection\thanks{%
Mathematics Subject Classifications: 05C38, 05C45.}}

\author{Bo Ning\thanks{%
Department of Applied Mathematics, Northwestern Polytechnical
University, Xi'an, Shaanxi 710072,
P. R. China}}
\maketitle

\begin{abstract}
In this short note we give a new sufficient condition for the existence of long cycles in graphs involving Fan-type degree conditions and neighborhood intersection.
\end{abstract}

\section{Introduction}

We use Bondy and Murty [2] for terminology and notation not defined here and consider simple graphs only.

Let $G$ be a graph, $H$ a subgraph of $G$, and $v$ a vertex of $G$. We use $N_H(v)$ to denote the set of neighbors of $v$ in $H$, and call $d_H(v)=|N_{H}(v)|$ the \emph{degree} of $v$ in $H$. For $x,y\in V(G)$, an $(x,y)$-\emph{path} is a path $P$ connecting $x$ and $y$, and vertices $x,y$ will be called the {\em end-vertices} of $P$. If $x,y\in V(H)$, the \emph{distance} between $x$ and $y$ in $H$, denoted by $d_H(x,y)$, is the length of a shortest $(x,y)$-path in $H$. When there is no danger of ambiguity, we will use $N(v)$, $d(v)$ and $d(x,y)$ instead of $N_G(v)$, $d_G(v)$ and $d_G(x,y)$, respectively.

Let $G$ be a graph and $G'$ be a subgraph of $G$. We call $G'$ an \emph{induced subgraph} of $G$ if $G'$ contains every edge $xy\in E(G)$ with $x,y\in V(G')$. A graph isomorphic to $K_{1,3}$ is called a \emph{claw}. A \emph{modified claw} is a graph isomorphic to one obtained by attaching an edge to one vertex of a triangle. We say that $G$ is \emph{claw-free} if $G$ contains no induced subgraph isomorphic to a claw.

In 1984, Fan [5] gave the following long cycle theorem involving the maximum degree of every pair of vertices at distance two in a 2-connected graph.

THEOREM 1 ([5]). Let $G$ be a 2-connected graph such that $\max{\{d(u),d(v)\}}\geq c/2$ for each pair of vertices $u$ and $v$ at distance 2. Then $G$ contains either a Hamilton cycle or a cycle of length at least $c$.

In [3], Bedrossian, Chen and Schelp gave an improvement of Fan's theorem. They further weakened the restriction on the pair of vertices $u$ and $v$ in graphs: they must be vertices of an induced claw or an induced modified claw at distance two.

THEOREM 2 ([3]) Let $G$ be a 2-connected graph such that $\max{\{d(u),d(v)\}}\geq c/2$ for each pair of nonadjacent vertices $u$ and $v$ that are vertices of an induced claw or an induced modified claw of $G$. Then $G$ contains either a Hamilton cycle or a cycle of length at least $c$.

On the other hand, Shi [7] gave the following sufficient condition for the existence of Hamilton cycles in claw-free graphs.

THEOREM 3 ([7]) Let $G$ be a 2-connected claw-free graph. If $|N(u)\cap N(v)|\geq 2$ for every pairs of vertices $u,v$ with $d(u,v)=2$, then $G$ contains a Hamilton cycle.

In this short note we give a new sufficient condition for the existence of long cycles involving Fan-type degree condition and neighborhood intersection. It can be seen as a generalization of Theorem 3.

THEOREM 4
Let $G$ be a 2-connected graph such that $\max\{d(u),d(v)\}\geq c/2$ for each pair of nonadjacent vertices $u$ and $v$ in an induced claw, and $|N(x)\cap N(y)|\geq 2$ for each pair of nonadjacent vertices $x$ and $y$ in an induced modified claw. Then $G$ contains either a Hamilton cycle or a cycle of length at least $c$.

\section{Proof of Theorem 4}
Before our proof, we give some additional useful terminology and notation. A \emph{$v_m$-path} is a path which has $v_m$ as one end-vertex. If a $v_m$-path is a longest path among all paths, then we call it a \emph{$v_m$-longest path}. Let $P=v_1v_2\ldots v_m$ be a path and denote by $t=t(P)=\max\{j:v_1v_j\in E(G)\}$.

The proof of Theorem 4 is motivated by [3]. It is mainly based on two lemmas below.

LEMMA 1 ([1]) Let $G$ be a 2-connected graph and $P$ be a longest path with two end-vertices $x$ and $y$. Then $G$ contains a Hamilton cycle or a cycle of length at least $d(x)+d(y)$.

LEMMA 2 Let G be a non-Hamiltonian 2-connected graph satisfying the condition of Theorem 4. Let $P=v_1v_2\ldots v_m (v_m=v)$ be a longest path of $G$. Then there exists a $v$-longest path such that the other end-vertex of the path has degree at
least $c/2$.

PROOF. Suppose not. Now we choose a path $P_1$ such that $t'=t(P_1)$ is as large as possible among all $v_m$-longest paths of $G$. Without loss of generality, we still denote $\overrightarrow{P_1}=v_1v_2\ldots v_m$.

CLAIM 1 $t'\leq m-1$.

PROOF. If $v_1v_m\in E(G)$, then $G$ is Hamiltonian or $G$ has a non-Hamilton cycle including all vertices of $P_1$. Hence $G$ has a Hamilton cycle or a path longer than $P_1$ since $G$ is 2-connected, a contradiction.

CLAIM 2. $\{v_1,v_{t'-1},v_{t'},v_{t'+1}\}$ induces a modified claw.

PROOF. By Claim 1 and  the choice of $t'$, $v_{t'+1}$ exists and $v_1v_{t'+1}\notin E(G)$. By the connectedness of $G$ and the choice of $P_1$, $t'\geq 3$. Assume that $v_1v_{t'-1},v_{t'-1}v_{t'+1}\notin E(G)$. Then $\{v_1,v_{t'-1},v_{t'},v_{t'+1}\}$ induces a claw. By the condition of Theorem 4 and the hypothesis that $d(v_1)<c/2$, we have $d(v_{t'-1})\geq c/2$ and $d(v_{t'+1})\geq c/2$. Let $P_{1}^{'}=v_{t'-1}\overleftarrow{P_1}v_1v_{t'}\overrightarrow{P_1}v_m$. Then
$P_1^{'}$ is a $v_m$-longest path such that the other end-vertex $v_{t'-1}$ with the degree at least $c/2$, a contradiction. If $v_{t'-1}v_{t'+1}\in E(G)$, then $P_1'=v_{t'-1}\overleftarrow{P_1}v_1v_t'\overrightarrow{P_1}v_m$ is a $v_m$-longest path with $t(P')\geq t'+1$, a contradiction.

By the assumption of Theorem 4, we have $|N(v_1)\cap N(v_{t'+1})|\geq 2$. By the definition of $t'$, there is a vertex $v_i\in N(v_1)\cap N(v_{t'+1})$, where $2\leq i \leq t'-2$.

CLAIM 3. $\{v_1,v_i,v_{i+1},v_{t'+1}\}$ induces a modified claw.

PROOF We have $v_1v_{i+1}\notin E(G)$, since otherwise $P_{1}^{'}=v_i\overleftarrow{P_1}v_{i+1}\overrightarrow{P_1}v_m$ is a $v_m$-longest path such that $t(P_{1}^{'})\geq t'+1$, a contradiction. If $v_{i+1}v_{t'+1}\notin E(G)$, then $\{v_1,v_i,v_{i+1},v_{t+1}\}$ induces an induced claw. By the condition of Theorem 4 and our hypothesis that $d(v_1)<c/2$, we have $d(v_{i+1})\geq c/2$. Let $P_{1}^{'}=v_{i+1}\overrightarrow{P_1}v_{t'}v_1\overrightarrow{P_1}v_iv_{t'+1}\overrightarrow{P_1}v_m$. Then $P_{1}^{'}$ is a $v_m$-longest path with the other end-vertex of degree at least $c/2$, a contradiction. Thus, we have $v_{i+1}v_{t'+1}\in E(G)$ and the proof of this claim is complete.

Now, we consider the same path $P_{1}^{'}=v_{i+1}\overrightarrow{P_1}v_{t'}v_1\overrightarrow{P_1}v_iv_{t'+1}\overrightarrow{P_1}v_m$. Then $P_{1}^{'}$
is a $v_m$-longest path such that $t(P_{1}^{'})\geq t'+1$, a contradiction. It is thus completed the proof of Lemma 2.

PROOF OF THEOREM 4. Suppose that $G$ contains no Hamilton cycles. By using lemma 2 twice, we obtain a longest path with both end-vertices having the degree at least $c/2$. Then by Lemma 2, we can find a cycle of length at least $c$.

\section{Concluding remarks}
In 1989, Zhu, Li and Deng [8] proposed the definition of implicit degree. Based on this definition, many theorems on long cycles including Theorems 1 and 2, are largely improved. For details, see [4,6,8]. Maybe it is interesting to find a version of Theorem 4 under the condition of implicit degree.

\section*{Acknowledgements}
This work is supported by NSFC (No.~11271300) and the Doctorate Foundation of Northwestern Polytechnical University (cx201326). The author thanks the anonymous referee for helpful comments on an earlier version of this note.

\end{document}